\documentclass[journal,twoside]{IEEEtran}

\usepackage{epstopdf}

\usepackage{cite}
\usepackage{amsmath,amssymb,amsfonts}
\usepackage{graphicx}
\usepackage{textcomp}
\usepackage{subcaption}
\usepackage[shortlabels]{enumitem}
\usepackage[ruled,vlined]{algorithm2e}

\usepackage{xspace}
\usepackage{soul}
\usepackage{mathtools}
\usepackage{makecell}
\usepackage{todonotes}
\usepackage{nicefrac}

\newcommand{\ra}{\rightarrow}

\newcommand{\ie}{\unskip, i.\,e.,\xspace}
\newcommand{\eg}{\unskip, e.\,g.,\xspace}

\newcommand{\sut}{\text{s.\,t.\,}}

\newcommand{\nrm}[1]{\left\lVert#1\right\rVert}

\newcommand{\abs}[1]{\left\lvert#1\right\rvert}

\newcommand{\tr}[1]{\ensuremath{\text{tr}}\left(#1\right)}
\newcommand{\E}[1]{\mathbb E\left[#1\right]}

\newcommand{\R}{\ensuremath{\mathbb{R}}}

\newcommand*\diff{\mathop{}\!\mathrm{d}}
\newcommand{\eps}{\ensuremath{\varepsilon}}

\newcommand{\spc}{\ensuremath{\,\,}}

\definecolor{dgreen}{rgb}{0.0, 0.5, 0.0}

\usepackage{float}
\usepackage{wrapfig}
\usepackage{xcolor}

\newtheorem{thm}{Theorem}

\renewcommand{\H}{\ensuremath{\mathcal{H}}}							

\makeatletter
\newcommand{\subalign}[1]{%
	\vcenter{%
		\Let@ \restore@math@cr \default@tag
		\baselineskip\fontdimen10 \scriptfont\tw@
		\advance\baselineskip\fontdimen12 \scriptfont\tw@
		\lineskip\thr@@\fontdimen8 \scriptfont\thr@@
		\lineskiplimit\lineskip
		\ialign{\hfil$\m@th\scriptstyle##$&$\m@th\scriptstyle{}##$\crcr
			#1\crcr
		}%
	}
}

\newcommand{\dt}{\ensuremath{\mathrm{d}t}}
\newcommand{\mt}{\ensuremath{^{\top}}}										
\newcommand{\inv}{\ensuremath{^{-1}}}										

\begin{document}

\title{A note on stabilizing reinforcement learning}
\author{Pavel Osinenko, Grigory Yaremenko, Ilya Osokin
}

\maketitle


\begin{abstract}
Reinforcement learning is a general methodology of adaptive optimal control that has attracted much attention in various fields ranging from video game industry to robot manipulators. 
Despite its remarkable performance demonstrations, plain reinforcement learning controllers do not guarantee stability which compromises their applicability in industry.
To provide such guarantees, measures have to be taken.
This gives rise to what could generally be called stabilizing reinforcement learning. 
Concrete approaches range from employment of human overseers to filter out unsafe actions to formally verified shields and fusion with classical stabilizing controllers. 
A line of attack that utilizes elements of adaptive control has become fairly popular in the recent years.
In this note, we critically address such an approach in a fairly general actor-critic setup for nonlinear time-continuous environments.
The actor network utilizes a so-called robustifying term that is supposed to compensate for the neural network errors.
The corresponding stability analysis is based on the value function itself.
We indicate a problem in such a stability analysis and provide a counterexample to the overall control scheme.
Implications for such a line of attack in stabilizing reinforcement learning are discussed.
Furthermore, unfortunately the said problem possess no fix without a substantial reconsideration of the whole approach.
As a positive message, we derive a stochastic critic neural network weight convergence analysis provided that the environment was stabilized.
\end{abstract}

\section{Introduction}

\subsection{Brief overview of reinforcement learning}
\label{sec:synopsis}


The design of a reinforcement learning agent is usually complex compared to the common techniques of nonlinear and adaptive control. 
A fairly standard basis for such a design is a time-discrete, frequently also discrete state- and action-space-based, environment \cite{Sutton2018}.
In this note, we consider a more general design that is suitable for time-continuous systems.
Namely, we will do this in the context of the work \cite{vamvoudakis2015asymptotically} that offers a stabilizing agent based on a robustifying term mentioned above.
So far, the details of the stabilization machinery are omitted -- they will be given later -- whereas now we are interested in the groundwork of the agent design.
We start with a continuous-time, deterministic, control-affine environment
\begin{equation}
\label{eqn:sys}
\dot x = F(x,u) := f(x) + g(x)u, x \in \R^n, u \in \R^m.
\end{equation}
An infinite-horizon optimal control problem\footnote{We should remark here that \cite{vamvoudakis2015asymptotically} addressed the case of input saturation, but this is non-essential for the current note} as follows is to be addressed:
\[
\min _{\kappa} J^\kappa\left(x_{0} \right)=\int_{0}^{\infty} \rho\left(x, \kappa(x)\right) \dt,
\]
where $\rho$ is a stage cost function (in other contexts -- reward, utility, instantaneous cost etc.).
The Hamiltonian of the problem reads, for a generic smooth function $h$: 
\begin{equation}
	\H(x, u | h): = \nabla h\mt(x) F(x,u) + \rho(x, u).
\end{equation}
Denoting the value function as $V := \min_\kappa J^\kappa$ and, accordingly, the optimal policy -- as $\kappa^*$, the Hamilton-Jacobi-Bellman equation (Hamilton-Jacobi-Bellman) reads:
\begin{equation}
\label{eqn:hjb}
\H\left(x, \kappa^{*}(x) | V\right)=0, \quad \forall x.
\end{equation}

Reinforcement learning, in a nutshell, amounts to finding approximate solutions $\hat V_\theta, \hat \kappa_\vartheta$ to \eqref{eqn:hjb}, as some $\theta$-, respectively, $\vartheta$-weighted neural networks, commonly called \textit{critic} and, respectively, \textit{actor}.
Let us assume one-hidden-layer topology for both.

Then, by an approximation theorem, one can express the value function and the optimal policy for some \textit{ideal} weights as follows: 
\begin{equation}
	V(x)={\theta^*}\mt \varphi(x)+\delta(x), \kappa^*(x) = {\vartheta^*}\mt \psi(x) + \delta_u(x)
\end{equation}
where $\varphi, \psi$ are the activation functions, and $\delta, \delta_u$ are the approximation errors.

The optimal policy has the property that the closed-loop of the agent and the environment (alias of the system \eqref{eqn:sys}), has an asymptotically stable origin.
Due to the presence of approximation errors $\delta(x),\delta_u(x)$, the actual agent policy $\hat \kappa_\vartheta$ cannot provide such a property in general.
Examples of divergence under reinforcement learning are known (see \eg \cite{fairbank2012divergence}) 
This leads to the need for special measures to remedy the issue and provide stability guarantees.
We briefly overview such measures in the next section.

\subsection{Brief overview of reinforcement learning}
Three general groups of approaches to stabilizing reinforcement learning can be formulated: shield-based reinforcement learning, fusion of model-predictive control and reinforcement learning and Lyapunov-based reinforcement learning.

The first group of approaches suggests to directly discard unsafe actions via, generally speaking, a filter which, depending on the context, may be called a shield, an overseer, a supervisor etc.
These algorithms vary in the way of how exactly the shield checks actions for safety and issues an emergency solution. 
The actual design of the shield may be as simple as human-based \ie manual \cite{saunders2017trial}, or as complex as formal-logic based \cite{Platzer2012-dyn-sys-log,Platzer2008-keymaera}.

Fusion with model-predictive control is active subfield in stabilizing reinforcement learning (see \eg \cite{zanon2019practical,Zanon2020,Koller2018,berkenkamp2019safe}).
It is no wonder since model-predictive control possesses such an established machinery for guaranteeing stability of the closed loop, such as terminal costs, terminal constraints, contraction tools etc.

Finally, the Lyapunov stability theory was recognized in reinforcement learning already by Perkins \& Barto \cite{perkins2001lyapunov}.
Advancements were made since \cite{chow2018,berkenkamp2017}, although such approaches to stabilizing reinforcement learning were offline and required checking a Lyapunov decay property over a region in the state space during running the agent.
This is not always feasible in industrial applications, such as robotics.

To this Lyapunov-based stabilizing reinforcement learning, agent design based on stabilization techniques of adaptive control can be categorized.
We focus in this note on the exemplary work \cite{vamvoudakis2015asymptotically}, although elements of such a methodology can be found in \eg \cite{zang2011robust,YangVamSafe}.
There is also a number of works inspired by \cite{vamvoudakis2015asymptotically} -- refer \eg to \cite{wei2020ImpulsiveADP, su2020decentralized, zang2020robustNonlinear, Jha2019LTIsystems, KinodynamicVamvoudakis, jha2018LTIwithoutPEcondition, MLbasedMultipleAgentsNguen}.
Essentially, a Lyapunov function for the closed-loop system was suggested in \cite{vamvoudakis2015asymptotically}, as well as one for the actor's and one for critic's weights, based on the value function utilizing a special \textit{robustifying term} in the agent's policy.
We give a synopsis of the work \cite{vamvoudakis2015asymptotically} in the next section.
Then, we indicate a flaw in the respective stability analysis.
It is stipulated that \cite{vamvoudakis2015asymptotically} erroneously followed the ideas of adaptive control by \cite{polycarpou2003line}.
Implications for similar control schemes are discussed.


%
%

\subsection{Actor-critic design}

To learn the actor and critic weights, one can follow what the Hamilton-Jacobi-Bellman prescribes.
From now on, let us assume that the action $u$ is generated by the policy $\hat \kappa_\vartheta$ (this is an \textit{on-policy} setup).
Then, the Hamiltonian under the critic reads: 
\begin{equation}
\H(x,u | \hat V_\theta) = \theta\mt \nabla \varphi(x) F(x, u) + \rho(x, u)
\end{equation}

Let us introduce the Hamiltonian approximation error as  
\begin{equation}
	\delta_\H(x, u) := \H(x,u | V) - \H(x,u | \hat V_{\theta^*}) = \nabla \delta\mt(x) F(x, u)
\end{equation}
By a similar token, consider the Hamiltonian \textit{temporal difference} (TD) as: 
\begin{multline}
	e_\H(\theta | x, u) := \H(x,u | \hat V_\theta) - \H(x,\kappa^* | V) =\\= \theta\mt \nabla \varphi(x) F(x,u) + \rho(x,u)
\end{multline}
Then, $\nabla_\theta e_\H (\theta | x,u) = \nabla \varphi(x) F(x,u)$.

Vamvoudakis et al. \cite{vamvoudakis2015asymptotically} suggested to use a \textit{data vector} 
\begin{equation}
	w(x,u) := \nabla \varphi(x) F(x,u)
\end{equation}
which yields 
\begin{equation}
e_\H(\theta | x, u) = \theta\mt w(x, u) + \rho(x, u)
\end{equation}
Introducing the weight error $\tilde \theta := \theta - \theta^*$, observe that $e_\H(\theta | x, u) = \tilde \theta\mt w(x, u) + \H(x, u | V) - \delta_H(x, u)$.
Let $\{x_{t_k}, u_{t_k}\}_{k=1}^M$ be an \textit{experience replay} of size $M$ and \sut $x_{t_M}, u_{t_M} = x(t), u(t)$ and $t_k < t_{k+1}, k \in [M]$, where $[M] = \{1, \dots, M\}$.
Normalizing the data vectors, one forms the critic loss as follows:
\begin{equation}
\label{eqn:critloss}
J_c(\theta | \{x_{t_k}, u_{t_k}\}_{k=1}^M ) := \frac 1 2 \sum_{k=1}^{M} \frac{e^2_\H(\theta | x_{t_k}, u_{t_k})}{ \left( w\mt_{t_k} w_{t_k} + 1 \right)^2 },
\end{equation}
where we used a shorthand notation $w_t : = w(x(t), u(t))$.
Then, the critic weights are learned in model-free fashion via online gradient descent as

\label{eqn:critupd}
\begin{multline}
\dot \theta :=  -\alpha \nabla_\theta J_c(\theta | \{x_{t_k}, u_{t_k}\}_{k=1}^M ) \\
 = - \alpha \sum_{k=1}^{M} \frac{ e_\H(\theta | x_{t_k}, u_{t_k}) \nabla_\theta e_\H(\theta | x_{t_k}, u_{t_k}) }{\left( w\mt_{t_k} w_{t_k} + 1 \right)^2} \\
 = - \alpha \sum_{k=1}^{M} \frac{ e_\H(\theta | x_{t_k}, u_{t_k}) w_{t_k} }{ \left( w\mt_{t_k} w_{t_k} + 1 \right)^2 }.
\end{multline}

Using the expression of the Hamiltonian TD through the Hamiltonian approximation error, one obtains:
\begin{equation}
\label{eqn:critupd_delta}
\dot \theta = -\alpha \sum_{k=1}^{M} \frac{w_{t_k} w_{t_k}\mt}{ \left( w\mt_{t_k} w_{t_k} + 1 \right)^2 } \tilde \theta + \alpha \sum_{k=1}^{M} \frac{w_{t_k} Z_{t_k} }{ \left( w\mt_{t_k} w_{t_k} + 1 \right)^2 },
\end{equation}
where $Z_{t_k} : = \H(x_{t_k}, u_{t_k} | V) - \delta_\H (x_{t_k}, u_{t_k})$.

Denote $\mathcal E_t : = \sum_{k=1}^{M} \frac{w_{t_k} w_{t_k}\mt}{ \left( w\mt_{t_k} w_{t_k} + 1 \right)^2 }$.
Then, one can perform a stability analysis of the differential equation \eqref{eqn:critupd_delta} under a suitable \textit{persistence of excitation} condition \ie 
\begin{equation}
	\mathcal E_t \succeq \eps I, \forall t \ge 0
\end{equation}
with $I$ being the identity matrix of proper dimension, and norm-boundedness of $Z_{t_k}$.
For the latter condition, it is crucial that the environment state $x$ be norm-bounded.
Regarding persistence of excitation, \cite{vamvoudakis2015asymptotically} suggested to update the experience replay until the respective condition is fulfilled.
In other words, one does not have to assume a kind of running persistence of excitation on the buffer.

Following an analogous methodology, one designs the actor to minimize an actor loss of the form 
\begin{equation}
	J_a( \vartheta | x ) = \frac 1 2 \tr{e_a(\vartheta | x)\mt e_a(\vartheta | x)}
\end{equation}
where the actor error is $e_a(\vartheta | x) := \hat \kappa_\vartheta(x) - \kappa^*(x) $.
The learning of the actor weights $\vartheta$ can also be done by stochastic gradient descent.

Concluding, the control scheme of \cite{vamvoudakis2015asymptotically} is fairly general, it is online \ie does not require pre-training of the agent, and it is plug-and-play in the sense that it does not require any special insight into the mathematical model of the environment \eqref{eqn:sys}.
However, it does not guarantee stability of the closed loop.
A special means is suggested to this end which we discuss in the next section.

\subsection{Issue with guaranteeing stability}
\label{sec:stability}

The full details of the actor can be found in \cite{vamvoudakis2015asymptotically}, whereas for us, the essential part is the \textit{robustifying term} that Vamvoudakis et al. used to guarantee closed-loop stability.
The ground idea comes from adaptive control \cite{polycarpou2003line}.
We depict it in a very short form.

Let us suppose that $x$ is scalar and $f, g$ are unknown, learned via models $\hat f(x) = \theta_f \varphi_f(x), \hat g(x) = \theta_g \varphi_g(x)$.
Let the ideal weights be $\theta^*_f, \theta^*_g$ and the respective weight errors be $\tilde \theta_f, \tilde \theta_g$. 
For brevity, we omit the $x$-argument from now on.
Let us consider a policy defined by $\kappa = \frac{1}{\hat g} (-K x - \hat f), K>0$.
Then, $\dot x = -K x + \tilde \theta_f \varphi_f + \tilde \theta_g \varphi_g \kappa$.

Taking weight update rules as $\dot \theta_f := \alpha_f x \varphi_f, \dot \theta_g := \alpha_g x \varphi_g \kappa$, and using a Lyapunov function candidate 
\begin{equation}
	L := \frac 1 2 x^2 + \frac 1 2 \tilde \theta_f\mt \alpha_f\inv \tilde \theta_f + \frac 1 2 \tilde \theta_g\mt \alpha_g\inv \tilde \theta_g
\end{equation}
one can show that $\dot L \le -K x^2$.
A projection might be required in the update of $\theta_g$ to ensure $\hat g$ be bounded away from zero, which is central to the described technique.

Now, \cite{vamvoudakis2015asymptotically} followed, roughly, this idea and added to the actor-generated policy $\hat \kappa_\vartheta$ a robustifying term of the form $- K \nrm{x}^2 \frac{I}{A + \nrm{x}^2}$.
The constants $K, A$ were then chosen to account for various bounds on the quantities involved (the activation function, its gradient, the Hamiltonian error etc.).
With this at hand, the value function $V$ was picked as the Lyapunov function candidate.
To be precise, the total Lyapunov function candidate in \cite{vamvoudakis2015asymptotically} consists of three parts: state-dependent, and two parts that depend on the actor and critic weight errors, respectively, but we are interested in the first one here.
Let $L$ denote this candidate.
Then, one has:
\begin{equation}
\label{eqn:LFVamvoudakis}
\begin{array}{lll}
\dot L & = & \nabla V\mt \left( f - g \tilde \vartheta\mt \psi + g(\kappa^* - \delta_u) - g  K \nrm{x}^2 \frac{I}{A + \nrm{x}^2} \right) \\
& = & -\rho^* - \nabla V\mt g \tilde \vartheta\mt \psi - \nabla V\mt g \delta_u - \\
& & \nabla V\mt g K \nrm{x}^2 \frac{I}{A + \nrm{x}^2},
\end{array}
\end{equation}
where $\rho^* := \rho(x, \kappa^*(x))$ and the last displayed identity follows from the Hamilton-Jacobi-Bellman $\nabla V\mt f = -\rho^* - \nabla V\mt g \kappa^*$.
Then, it was erroneously claimed that\footnote{See equation (45) in \cite{vamvoudakis2015asymptotically}}
\begin{equation}
\label{eqn:LFVamvoudakis_fail}
\begin{array}{l}
-\rho^* - \nabla V\mt g \tilde \vartheta\mt \psi - \nabla V\mt g \delta_u -\nabla V\mt g K \nrm{x}^2 \frac{I}{A + \nrm{x}^2} \le \\
-\rho^* - \overline {\nabla V} \bar g \bar \psi \nrm{\tilde \vartheta} - \overline {\nabla V} \bar g \bar \delta_u - \overline {\nabla V} \bar g K \nrm{x}^2 \frac{I}{A + \nrm{x}^2},
\end{array}
\end{equation}
where the lines over variables indicate bounds on the respective norms or absolute values.
Whereas it might have been a typo in the second and third summand in the last displayed expression (they should come with a plus sign), the central \textit{issue} is with the last term.
It is exactly this term that should have eventually served the closed-loop stability.
But, one cannot claim that it comes with a minus sign since $\nabla V\mt g K \nrm{x}^2 \frac{I}{A + \nrm{x}^2}$ is sign-indefinite, not to mention that $g$ may actually vanish (cf. boundedness away from zero in the above described adaptive control case).

\section{Discussion and implications}
\label{sec:discussion}

The most severe ramification of the above issue is that the designed reinforcement learning agent evidently does not provide stability guarantees.
One might argue that the source of the problem is the use of the value function as a Lyapunov function candidate.
Generally speaking, one obtains in this case an expression similar to $\dot L = -\rho^* - \nabla V g \kappa^* + \nabla V g \kappa$.
One can see that no matter how $\kappa$ is chosen, the last term is, in general, sign-indefinite.
One could make use of a bound $\nrm{\kappa - \kappa^*}$ to benefit from at least a portion of the decay via $-\rho^*$, but \textit{the fundamental issue with online reinforcement learning is that one cannot know in advance how good the chosen neural network topology is capable of approximating the optimal policy}.
Therefore, such a bound is a priori not available in general.

There is another, conceptual, issue with the robustifying term as described above, namely, its interference with the nominal actor.
The gain $K$ might need to be chosen large enough depending on the various assumed bounds (if the overall routine had been correct).
Such a setting might have turned out conservative.
But this could have harmed a key principle of stabilizing or, more generally speaking \textit{safe} reinforcement learning -- \textit{minimal interference} -- that says that stability and/or safety measures (like a robustifying term in the considered case) should not interfere with the actor's learning.
Thus, even if one had followed the methodology of stabilization via robustifying controls correctly, one would have also had to account for possible influence of it onto the agent's learning.
Concluding this section, unfortunately, no trivial fix to the stability guarantee issue of the discussed reinforcement learning agent based on a robustifying term exists.
In the next section, we give a simple counterexample to the methodology suggested in \cite{vamvoudakis2015asymptotically}.

\subsection{A counterexample}

In this section we will follow the work \cite{diffGamesDoyle} and the principles of converse of optimality.
Let us consider the following system 
\begin{equation}
\begin{cases}
	\dot{x}_{1}=x_{2} \\
	\dot{x}_{2}=f(x)+g(x) u,\\
\end{cases}
\end{equation}
for which we assume a stage cost $\rho (x)=x_{2}^{2} + u^2$ and a Lyapunov function candidate $L$ equal to the value function \ie $L:=V = x_1^2 + x_2^2$.
The directional derivative of $L$ along $F$ (in the sense of equation \ref{eqn:sys}) reads:
\begin{equation}
\mathcal{L}_{F}L = \langle\nabla L, F\rangle= x_{2} \frac{\partial V}{\partial x_{1}}+(f(x)+g(x) u) \frac{\partial V}{\partial x_{2}}.
\end{equation}
Then, 
\begin{multline}
\nabla V=2\left(\begin{array}{c}
x_{1} \\
x_{2}
\end{array}\right)\Rightarrow \\ \mathcal{L}_{F} L=2 x_{1} x_{2}+2(f(x)
+g(x) u) x_{2}.
\end{multline}
The Hamilton-Jacobi-Bellman equation reads: 
\begin{equation}
	V_{1} x_{2}+V_{2} f(x)-\frac{1}{4}\left(V_{2}\right)^{2} g^{2}(x)+\rho (x)=0.
\end{equation}
where $V_{i} :=\frac{\partial V(x)}{\partial x_{i}}$.
This can be cast into
\begin{equation}\label{eqn:Hamilton-Jacobi-Bellman}
2 x_{1} x_{2}+2 x_{2} f(x)=\frac{1}{4}\left(2 x_{2}\right)^{2} g^{2}(x)-x_{2}^{2}.
\end{equation}
From \eqref{eqn:Hamilton-Jacobi-Bellman}, we can find that 
\begin{equation}
	f(x)=-x_{1}-\frac{1}{2} x_{2}\left(1-g^{2}(x)\right).
\end{equation}
Therefore, we may rewrite the original system as
\begin{equation}\left\{\begin{array}{l}
	\dot{x}_{1}=x_{2} \\
	\dot{x}_{2}=-x_{1}-\frac{1}{2} x_{2}\left(1-g^{2}(x)\right)+g(x) u
	\end{array}\right.
\end{equation}
The optimal policy reads $\kappa^{*}(x)=-g(x) x_{2}$.
Let us assume that the actor somehow learned this policy exactly and we had a robustifying term.
The effective policy applied then reads: 
\begin{equation}
	\kappa(x) = \kappa^*(x) - \frac{K\|x\|^{2}}{A+\|x\|^{2}}  = - g(x) x_{2} - \frac{K\|x\|^{2}}{A+\|x\|^{2}},
\end{equation}
where $ - \frac{K\|x\|^{2}}{A+\|x\|^{2}}$ is the robustifying term.
Substituting the effective policy into $\mathcal{L}_{F}L$ yields:
	\begin{multline}
		\mathcal{L}_{F} L = -x_2^2 - x_2^2 g^2(x) \underbrace{- 2 x_2 g(x) \frac{K \nrm{x}^2}{A + \nrm{x}^2}}_{=:R}.		
		\end{multline}
Let us now analyze what is going on here.
First off, in absence of the robustifying term, the derivative would equal just the stage cost.
The term $R$ is induced by the robustifying effect.
Suppose, for the sake of demonstration, that $g(x) \equiv 1, x_1=0, A=1$.
Then,
\[
	R = - 2 x_2 \frac{K x_2^2}{A + x_2^2}.
\]
Under the assumption that $g(x) \equiv 1, x_1=0, A=1$, it can then be deduced that $\mathcal{L}_{F} L$ is positive in the set
\[
	x_2 \in \left[ \frac 1 2 (-K - \sqrt{K^2 - 4}), 0 \right)
\]
which gets larger with larger robustifying gain $K$.
In the said set, no decay of the Lyapunov function candidate can be guaranteed.			

Thus, we have shown that the approach, proposed by Vamvoudakis et al. does not stabilize the system with guarantees.
In contrary, in a simple system, even in the case when the actor has learned the optimal policy perfectly, the effective policy with a robustifying term acts in a destabilizing way.

Now, we proceed to the discussion of the critic part which is also flawed, but less severely.

\section{Issue with the generic-action Hamiltonian}
\label{sec:hamiltonian}

There is another subtlety in \cite{vamvoudakis2015asymptotically} that has to do with the critic learning.
This issue is two-fold: first, exponential convergence of the critic weight error (see equation 25 in the cited work) was claimed, while it should have been ultimate boundedness instead, taking into account the perturbation term (the second one in equation \ref{eqn:critupd_delta} above).
Second, the bound on the perturbation term was incorrect.
Let us recall \eqref{eqn:critupd_delta}.
The perturbation term here reads:
\begin{equation}
\label{eqn:perturb}
\alpha \sum_{k=1}^{M} \frac{w_{t_k} \left( \H(x_{t_k}, u_{t_k}|V) - \delta_\H (x_{t_k}, u_{t_k}) \right) }{ \left( w\mt_{t_k} w_{t_k} + 1 \right)^2 }.
\end{equation}
In \cite{vamvoudakis2015asymptotically}, the generic-action Hamiltonian $\H(x_{t_k}, u_{t_k}|V)$ was erroneously assumed zero, while it is only zero under the optimal policy (see equation \ref{eqn:hjb}).
What this \textit{implies} is that the critic learning quality actually depends on the policy applied to the environment.
In the case of \cite{vamvoudakis2015asymptotically}, the overall routine is an on-policy reinforcement learning, which means that the policy being learned is the one applied.

Ideally, one would expect the respective generic-action Hamiltonian approach zero as the learned policy approaches the optimal one.
Otherwise, the ultimate critic weight error will depend on the quality of the actual policy.
Notice that application of robustifying terms, in general, worsens this quality.
Furthermore, of course, the critic's weight error convergence collapses when the environment is unstable.
In the next section, we fix the discussed issue, extend the analysis to the case of a stochastic environment, while \textit{assuming} closed-loop stability (disregarding how a stabilizing policy was achieved).

\section{New result: stochastic critic convergence}

In this section, we study a similar actor-critic design as before, but for an environment described by a stochastic differential equation 
\begin{equation}
	\diff X_{t} = f\left(X_{t}, U_{t}\right) \dt + \sigma\left(X_{t}, U_{t}\right) \mathrm d B_{t}
\end{equation}
where $\{B_t\}_{t>0}$ is a vector Brownian motion.

The optimal control problem is now in the following format:
\begin{multline}
\label{eqn:ocp_stoch}
\min_{\kappa} J^\kappa(x_0)\\
 = \E{ \int_{0}^{\infty} e^{-\gamma t} \rho\left(X_{t}, \kappa\left(X_{t}\right)\right) \dt \mid X_{0}=x_{0} },
\end{multline}
where we introduced a discounting factor $\gamma$ to be more general.
The Hamiltonian now reads, for a generic smooth function $h$,
\begin{align*}
& \H(x,u |h) = \nabla h(x)\mt f(x,u) \\
&  \qquad + \frac 1 2 \tr{\sigma\mt(x,u) \nabla^2 h(x) \sigma(x,u)} + \rho(x,u) - \gamma h(x),
\end{align*}
while the Hamilton-Jacobi-Bellman is the same as above.
Now, the Hamiltonian under the critic reads:
\begin{align*}
&\H(x,u | \hat V_\theta) = \\
&\theta\mt \nabla \varphi(x) f(x, u) + \frac  1 2 \theta\mt \eta(x,u) +  \rho(x,u) - \gamma \theta\mt \varphi(x),
\end{align*}
where $\eta(x,u)$ is the vector with entries $\tr{\sigma\mt(x,u) \nabla^2 \varphi_i(x) \sigma(x,u) }$, $i \in [N_c]$, where $N_c$ is the number of the critic's features and the notation $\varphi_i(x)$ means the $i$th feature.
The Hamiltonian approximation error in our case amounts to
\begin{align*}
& \delta_\H(x,u) = \\
& \nabla \delta(x)\mt f(x,u) + \frac 1 2 \tr{\sigma\mt(x,u) \nabla^2 \delta(x) \sigma(x,u)} - \gamma \delta(x).
\end{align*}
The Hamiltonian TD modifies to 
\begin{multline}
	e_\H(\theta | x,u) = \theta\mt \nabla \varphi(x) f(x,u) \\+ \frac 1 2 \theta\mt \eta(x,u) + \rho(x,u) - \gamma \theta\mt \varphi(x)
\end{multline}
The data vector in turn becomes 
\begin{equation}
w(x,u) := \nabla \varphi(x) f(x, u) + \frac 1 2 \eta(x, u) - \gamma \varphi(x)
\end{equation}
Then, $e_\H(\theta | x,u) = \tilde \theta\mt w(x, u) + \H(x,u|V) - \delta_\H(x,u)$.

Following the same procedure as in the first section of this note, the evolution of the critic weight error in norm square follows the following SDE:
\begin{equation}
\label{eqn:critupd_stoch}
\begin{array}{ll}
\mathrm d \nrm{\tilde \Theta_t}^2 = & -\alpha \bigg( \nrm{\tilde \Theta_t}^2_{\mathcal E_t} + \\
& \sum \limits_{k=1}^{M} \frac{\tilde \Theta_t W_{t_k} \left(\H(X_{t_k}, U_{t_k} | V) - \delta_\H (X_{t_k}, U_{t_k})\right) }{ \left( W\mt_{t_k} W_{t_k} + 1 \right)^2 } \bigg) \dt, 
\end{array}
\end{equation}
where $\nrm{\bullet}_{\bullet}$ means the weighted Euclidean norm and 
\begin{equation}
	\mathcal E_t := \sum_{k=1}^{M} \frac{W_{t_k} W_{t_k}\mt}{ \left( W\mt_{t_k} W_{t_k} + 1 \right)^2 }
\end{equation}


\begin{thm}
Consider the critic learning via \eqref{eqn:critupd_stoch} under some behavior policy $\mu$.
Let the following conditions hold:
\begin{enumerate}[(a)]
\item \label{a} $f(\bullet, \mu(\bullet))$ and $\sigma(\bullet, \mu(\bullet))$ are of linear growth and $\mu(0) = 0, f(0,0) = 0, \sigma(0,0) = 0$;
\item \label{b} the critic topology is chosen \sut the value function approximation error $\delta$ is of quadratic growth and $\delta(0) = 0$;
\item \label{c} $\rho$ and $V$ are of quadratic growth;
\item \label{d} the behavior policy is persistently exciting \ie $\mathcal E_t \succeq \eps I$ a.~s.;  
\item \label{e}$\exists \bar X > 0 \spc \forall t \ge 0 \spc \E{\nrm{X_t}^4} \le \bar X^2$.
\end{enumerate}
Then, the critic error weights satisfy:
\begin{equation}
\label{eqn:critwgt_bound_final}
\mathbb{E}\left[\left\|\widetilde{\Theta}_{t}\right\|^{2}\right] \leqslant e^{-\alpha \varepsilon t}\|\tilde{\theta}(0)\|^{2}+D \sup _{\tau \leq t} \sqrt{\mathbb{E}\left[\left\|\tilde{\Theta}_{\tau}\right\|^{2}\right]} 
\end{equation}
for some constant $D$.
\end{thm}

\begin{IEEEproof}
Consider \eqref{eqn:critupd_stoch}.
Denote $Z_{t_k} := \H(X_{t_k}, U_{t_k} | V) - \delta_\H (X_{t_k}, U_{t_k})$.
Now, introduce a stopping time $T_R := \inf\limits_{t>0}\bigg\{\nrm{X_t}>R\bigg\}$ and integrate \eqref{eqn:critupd_stoch} from $0$ to $\mathcal T_R := t \land T_R$, where $t \land T_R := \min\{t, T_R\}$ and take expectation on both sides of \eqref{eqn:critupd_stoch}:
\begin{equation}
\label{eqn:critupd_stoch_stp_time}
\begin{array}{ll}
\E{ \int\limits_0^{\mathcal T_R} \nrm{ \tilde \Theta_\tau}^2 \diff \tau } = \\
\quad \E{ \int\limits_0^{\mathcal T_R}  -\alpha \left( \nrm{\tilde \Theta_t}^2_{\mathcal E_t} + 
 \sum \limits_{k=1}^{M} \frac{\tilde \Theta_t W_{t_k}Z_{t_k} }{ ( W\mt_{t_k} W_{t_k} + 1)^2 } \right) \diff \tau} \le \\
e^{-\alpha \eps \mathcal T_R} \nrm{\tilde \theta(0)} + \frac \alpha 2 e^{ - \alpha \mathcal T_R} \E{ \int_{0}^{\mathcal T_R} e^{\alpha \tau}\left\|\tilde{\Theta}_{\tau}\right\| \sum \limits_{k=1}^{M} Z_{t_k} \mathrm d \tau},
\end{array}
\end{equation}
where the last inequality follows from (d) the fact that $\frac{\nrm{w}}{(w\mt w +1)^2} \le \frac 1 2$ always.
Now, due to (e), $\mathcal T_R \ra t$ a.~s. as $R \ra \infty$, and so using the dominated convergence on the right of \eqref{eqn:critupd_stoch_stp_time} and the Fatou's lemma on the left of \eqref{eqn:critupd_stoch_stp_time}, deduce:
\begin{equation}
\label{eqn:critwgt_bound}
\begin{array}{ll}
\E{\nrm{\tilde \Theta_t}^2} \le & e^{-\alpha \eps t} \nrm{\tilde \theta(0)} \\
& + \frac \alpha 2 e^{- \alpha t} \E{ \int_{0}^{t} e^{\alpha \tau}\left\|\tilde{\Theta}_{\tau}\right\| \sum \limits_{k=1}^{M} Z_{t_k} \mathrm d \tau}.
\end{array}
\end{equation} 
Specifically, the respective intermediate steps read:
\[
\begin{aligned}
\E{\nrm{\tilde \Theta_t}^2} \le \liminf_{R \ra \infty} \E{ \int\limits_0^{\mathcal T_R} \nrm{ \tilde \Theta_\tau}^2 \diff \tau } \\ \text{(Fatou's lemma)}
\end{aligned}
\]
\[
\begin{aligned}
\lim_{R \ra \infty} \E{ \int_{0}^{\mathcal T_R} e^{\alpha \tau}\left\|\tilde{\Theta}_{\tau}\right\| \sum \limits_{k=1}^{M} Z_{t_k} \mathrm d \tau} = \\
 \E{ \lim_{R \ra \infty} \int_{0}^{\mathcal T_R} e^{\alpha \tau}\left\|\tilde{\Theta}_{\tau}\right\| \sum \limits_{k=1}^{M} Z_{t_k} \mathrm d \tau } \\ \text{(dominated convergence)}
\end{aligned}
\]
Combining these facts with \eqref{eqn:critupd_stoch_stp_time} gives \eqref{eqn:critwgt_bound}.

From \ref{a}, \ref{b}, \ref{c} we can see that $\exists C>0$ \sut $\abs{\H(x,\mu(x)|V) - \delta_\H(x, \mu(x))} \le C \nrm{x}^2$.

From \ref{e}, $\forall t \ge 0 \spc \E{ \abs{Z_t}^2 } \le C^2 \bar X^2$.
Then, using the Fubini's lemma and the Cauchy-Schwartz inequality, one obtains from \eqref{eqn:critwgt_bound}:
\begin{multline}
\label{eqn:critwgt_bound_final_prf}
\mathbb{E}\left[\left\|\widetilde{\Theta}_{t}\right\|^{2}\right] \leqslant e^{-\alpha \varepsilon t}\|\tilde{\theta}(0)\|^{2}\\+\frac{M C \bar X}{2} \sup _{\tau \leq t} \sqrt{\mathbb{E}\left[\left\|\tilde{\Theta}_{\tau}\right\|^{2}\right]},	
\end{multline}
where $\frac{M C \bar X}{2}$ is the desired constant $D$ from the theorem statement. 
\end{IEEEproof}
From \eqref{eqn:critwgt_bound_final_prf} one can deduce a suitable ultimate boundedness property for the mean-square critic weight error.
A remark should be made: the above stated assumptions are only necessary when we are dealing with an SDE driven by Brownian motion.
If the driving noise were bounded, the assumptions could be discarded.
\section{Conclusion}
\label{sec:conclusion}

Guaranteeing system stability under reinforcement learning remains a challenge.
Various classical control techniques are frequently invoked to this end.
There is a fairly general methodology that uses robustifying controls, which in turn does not directly utilize the environment model, to compensate for the agent's neural network approximation errors.
This note discussed issues with the corresponding stability analysis of the agent-environment closed loop.
Stabilizing reinforcement learning agent designs that base on robustifying controls in the similar manner, as discussed above, should thus be carried out with precautions.


\bibliographystyle{IEEEtran}
\bibliography{
bib/MPC,
bib/stability,
bib/RL,
bib/robotics,
bib/Osinenko,
bib/formal,
bib/stochastic,
bib/nonsmooth,
bib/adaptive,
bib/bad
}



\end{document}